\documentclass[11pt,draft]{amsart}
\usepackage{dsfont}
\usepackage{amsfonts}
\usepackage{mathrsfs}
\usepackage{amssymb}
\usepackage{amsmath}
\usepackage{amsfonts}
\usepackage{amsfonts,enumerate}
\usepackage{pifont}
\usepackage{enumerate}
\usepackage{graphics}
\usepackage{verbatim}
\usepackage{amsthm}
\usepackage{latexsym,bm}
\usepackage{color}
\numberwithin{equation}{section}
\newtheorem{theorem}{Theorem}[section]

\theoremstyle{definition}
\newtheorem{remark}[theorem]{Remark}

\newcommand{\R}{{\mathbb R}}
\newcommand{\Z}{\mathbb Z}

\newcommand{\E}{{\mathcal E}}

\newcommand{\un}{{\mathds {1}}}

\newcommand{\I}{{\mathcal{I}}}

\begin{document}

\title[]{The behaviour of square functions from ergodic theory in $L^{\infty}$}

\author{Guixiang Hong}
\address{Instituto de Ciencias Matem\'aticas,
CSIC-UAM-UC3M-UCM, Consejo Superior de Investigaciones
Cient\'ificas, C/Nicol\'as Cabrera 13-15. \newline 28049, Madrid. Spain.\\
\emph{E-mail address: guixiang.hong@icmat.es}}
\thanks{The author is supported by MINECO: ICMAT Severo Ochoa project SEV-2011-0087 and ERC Grant StG-256997-CZOSQP (EU)}

\thanks{\small {{\it MR(2000) Subject Classification}.} Primary
42B25; Secondary 47G10.}

\maketitle

\begin{abstract}
In this paper, we analyze carefully the behaviour in $L^\infty(\R)$ of the square functions $S$ and $S_\I$'s, originating from ergodic theory. Firstly, we show that we can find some function $f\in L^\infty(\mathbb{R})$ , such that $Sf$ equals infinity on a nonzero measure set. Secondly, we can find compact supported function $f\in L^\infty(\mathbb{R})$ and $\I$ such that $S_\I f$ does not belong to $BMO$ space. Finally, we show that $S$ is bounded from $L^{\infty}_c$ to $BMO$ space. As a consequence, we solve an open question posed by Jones, Kaufman, Rosenblatt and Wierdl in \cite{JKRW98}. That is, $S_\I$ are uniformly bounded in $L^p(\R)$ with respect to $\I$ for $2<p<\infty$.
\end{abstract}

\section{Introduction}
A variety of square functions were introduced in \cite{JKRW98} by Jones {\it et al} as tools to deal with variational inequalities, whence measure the speed of the convergence of a sequence of differential averages. To present the square functions we are interested in this paper, we need some notations. Let $\sigma_k$ be the $k$-th dyadic $\sigma$ algebra in $\R$. That is, $\sigma_k$ is generated by the dyadic intervals with side-length equal to $2^k$. Denote by $\E_k$ the expectation with respect to $\sigma_k$. For $x\in\R$,
let $I_k(x)$ denote any possible interval containing $x$ with length $2^k$.
Let $\mathcal{I}=\{I_k(x)\}_{k\in\mathbb{Z},x\in\R}$, for any finite compact supported function $f$ on $\R$, define
\begin{align}\label{def of si}
S_\I f(x)=\big(\sum_{k\in\Z}|M_{I_k(x)}f(x)-\E_kf(x)|^2\big)^{1/2},
\end{align}
where
$$M_{I_k(x)}f(x)=2^{-k}\int_{I_{k}(x)}f(y)dy.$$

In Theorem 2.2 of \cite{JKRW98}, the authors proved that $S_\I$ is bounded in $L^2(\mathbb{R})$ uniformly with respect to $\I$. That is, there exist a constant $C>0$ independent of $\I$ such that
\begin{align}\label{l2i}
\|S_\I f\|_{2}\leq C\|f\|_{2},\;\forall f\in L^{2}(\mathbb{R}).
\end{align}
But in Remark 4.5 of the same paper, the authors observed that for some $\I$, $S_\I$ may not map $L^{\infty}$ to $BMO_d$ (the dyadic BMO space on the torus). Hence the interpolation argument can not be applied, and they leave it as an open question (Question 4.7 in the same paper) that whether $S_\I$ is bounded in $L^p(\mathbb{R})$ uniformly with respect to $\I$ for $2<p<\infty$.

In this paper, we give a positive answer of this question. In order to present our approach, we need more notations. Let $\I_k$ be the set of intervals containing the origin with length $2^k$. We will consider the following square function
\begin{align}\label{def of s}
Sf(x)=\big(\sum_{k\in\Z}\sup_{I\in\I_k}|M_{I+x}f(x)-\E_kf(x)|^2\big)^{1/2}.
\end{align}
It is clear that for all $\I$, $S_\I f(x)\leq Sf(x)$ for all $f$ and almost every $x\in\R$. Hence for all $1<p<\infty$, $L^p$-boundedness of $S$ implies the uniform $L^p$-boundedness of $S_\I$, since the spaces $L^p$'s are K\"othe function spaces. Hence it suffices to prove $S$ is bounded on $L^p(\mathbb{R})$ for all $2<p<\infty$.

It is known from Theorem A' of \cite{JRW03} that $S$ is bounded on $L^2(\mathbb{R})$, i.e.
\begin{align}\label{l2}
\|Sf\|_{2}\leq C\|f\|_{2},\;\forall f\in L^{2}(\mathbb{R}),
\end{align}
for some positive constant $C$. Hence the $L^p$-boundedness would be obtained by interpolation, if we could show the $(L^\infty,BMO_d)$-boundedness of $S$. However, as recalled previously, for some $\I$, $S_{\I}$ may not be bounded from $L^\infty$ to $BMO_d$. Now a key observation is that the $(L^\infty,BMO_d)$-boundedness of $S$ may still survive, since $BMO$ spaces are not K\"other function spaces.

Therefore, in Section 2, we carefully analyze the behavior of $S_\I$ and $S$ in $L^{\infty}$. The results we obtain can be concluded as follows.
\begin{theorem}\label{main theorem 1}
{\rm (i)} There exist a function $f\in L^{\infty}(\R)$ and a nonzero measure set $E\subset\R$ such that for any $x\in E$, $S_\I f(x)=\infty$ for some $\I$, whence $Sf(x)=\infty$.

{\rm (ii)} There exist a compact supported function $f\in L^{\infty}(\R)$ and $\I$ such that $S_\I f\notin BMO_d(\R)$.
\end{theorem}

The {\rm (ii)} is interesting in the sense that it is different from\cite{Wan85}, where the author proved that the classical $g$-function is bounded from $L^{\infty}_c(\R)$, the space of compact supported $L^\infty(\R)$ functions, to $BMO(\R)$ even thougth we can find $f\in L^{\infty}(\mathbb{R})$ such that $g(f)=\infty$ almost everywhere.

On the other hand,
From {\rm (i)}, we can not expect $S$ maps the whole $L^\infty$ to $BMO_d$. On the other hand, by the $L^2$-boundedness (\ref{l2}), for almost every $x\in\R$, $Sf(x)<\infty$ for all $f\in L^\infty_c(\R)$, . Hence the best result we can expect is the following, which will be shown in Section 3.

\begin{theorem}\label{main theorem 2}
$S$ is bounded from $L^{\infty}_c(\R)$ to $BMO_d(\R)$. Hence by interpolation, for all $2<p<\infty$, $S$ is bounded on $L^p(\R)$, whence $S_\I$ is uniformly bounded on $L^p(\R)$.
\end{theorem}

\section{Proof of Theorem \ref{main theorem 1}}
This section is devoted to the proof of Theorem \ref{main theorem 1}. Some comments on the cases $\mathbb{Z}$, $\mathbb{T}$ and $\mathbb{R}^n$ with $n\geq2$ are also included at the end of this section.
\begin{proof}
The proof of {\rm (i)}. Take $f=\chi_{[0,\infty)}$ and $I_k(x)=[x-2^k,x)$ for any $k\in\Z$. For any $x\in[0,\infty)$, there exist $\ell>0$ such that $x\in[0,2^{\ell})$.

Fix a $x\in[0,2^{\ell})$. Obviously, $\E_kf(x)=1$, for all $k\in\Z$. On the other hand, if $k>\ell$, then $x-2^k<0$. Thus
$$M_{I_k(x)}f(x)=2^{-k}\int^x_{x-2^k}f=2^{-k}x\leq {2}^{-1}.$$
Therefore,
\begin{align*}
S_{\I}f(x)&\geq\big(\sum_{k>\ell}|M_{I_k(x)}f(x)-\E_kf|^2\big)^{1/2}\\
&\geq\big(\sum_{k>\ell}|1-1/2|^2\big)^{1/2}=\infty.
\end{align*}

The proof of (ii). The basic construction is similar to that in Remark 4.5 of \cite{JKRW98}, where the authors proved in the torus case. Let $I_{\ell}=[1/2,1/2+1/2^{\ell})$. Let $P$ and $N$ denote two disjoint subsets of $[1/2,1)$ such that for all $\ell>2$, $$|P\cap I_{\ell}|=|N\cap I_\ell|=2^{-(\ell+1)}.$$
Take $I_k(x)=(x-2^{k},x]$ for each $x\in N$, and  $I_k(x)=(x,x+2^{k}]$ for each $x\in P$. Let $f=\chi_{[1/2,1)}$ and let $\ell>2$.

Fix $x\in I_\ell$. $\E_kf(x)=1$ for $k\leq-2$; $\E_kf(x)=2^{-k-1}$ for $k\geq-1$. Moreover, if $x\in P$, then
$$M_{I_k(x)}f(x)=2^{-k}\int^{x+2^k}_xf=1,\;\mathrm{for}\;k\leq-2;$$
$M_{I_k(x)}f(x)=2^{-k}(1-x)$ for $k\geq-1$.
If $x\in N$, then for $k>-\ell$, $x-2^k< 1/2$. Thus
$$M_{I_k(x)}f(x)=2^{-k}\int^x_{x-2^k}f=2^{-k}(x-1/2)\leq 2^{-k-\ell}<1/2.$$

To conclude, for $x\in P\cap I_\ell$, we have
\begin{align*}
S_{\I}f(x)&\leq\big(\sum_{k\in\Z}|M_{I_k(x)}f(x)-\E_kf|^2\big)^{1/2}\\
&=\big(\sum_{k\geq-1}|M_{I_k(x)}f(x)-\E_kf|^2\big)^{1/2}\\
&=\big(\sum_{k\geq-1}|2^{-k}(x-1/2)|^2\big)^{1/2}\leq1.
\end{align*}
While for $x\in N\cap I_\ell$, we have
\begin{align*}
S_{\I}f(x)&\geq\big(\sum_{\ell<k\leq-2}|M_{I_k(x)}f(x)-\E_kf|^2\big)^{1/2}\\
&=\big(\sum_{\ell<k\leq-2}|1-1/2|^2\big)^{1/2}= \frac{1}{2}\sqrt{\ell-1}.
\end{align*}
Then, it is easy to check that for for large $\ell$,
$$\big|\int_{N\cap I_\ell}(Sf(x)-Sf(y))dy\big|\geq2\big|\int_{P\cap I_\ell}(Sf(x)-Sf(y))dy\big|$$
for any $x\in P\cap I_\ell$.
Therefore, by triangle inequalities
\begin{align*}
&\|Sf\|_{BMO_d}\geq \frac{1}{|I_\ell|^2}\int_{I_\ell}\big|\int_{I_\ell}(Sf(x)-Sf(y))dy\big|dx\\
&\geq \frac{1}{|I_\ell|^2}\int_{P\cap I_\ell}\big|\int_{N\cap I_\ell}(Sf(x)-Sf(y))dy+\int_{P\cap I_\ell}(Sf(x)-Sf(y))dy\big|dx\\
&\geq \frac{1}{|I_\ell|^2}\int_{P\cap I_\ell}\Big(\big|\int_{N\cap I_\ell}(Sf(x)-Sf(y))dy\big|-\big|\int_{P\cap I_\ell}(Sf(x)-Sf(y))dy\big|\Big)dx\\
&\geq \frac{1}{2}\frac{1}{|I_\ell|^2}\int_{P\cap I_\ell}\big|\int_{N\cap I_\ell}(Sf(x)-Sf(y))dy\big|dx\\
&\geq \frac{1}{2}\frac{1}{|I_\ell|^2}\int_{P\cap I_\ell}\int_{N\cap I_\ell}(\frac{1}{2}\sqrt{\ell-1}-1)dydx\geq\frac{1}{16}\sqrt{\ell}.
\end{align*}
This finishes the proof since $\ell$ can be taken as large as we want.
\end{proof}

In the case $\mathbb{T}$, $\Z$ and $\R^n$, we can define $S_\I$ (or $S_{\mathcal{Q}}$) and $S$ similarly.

\begin{remark}
{\rm (i)}. The case $\Z$. Take $f=\chi_{[0,\infty)}$ and $I_k(j)=[j-2^k,j)$ for any $k\in\mathbb{N}$, then using the same arguments, we can show an analog of Theorem \ref{main theorem 1} (i). However Theorem \ref{main theorem 1} (ii) never be true in this case. Actually it follows from Theorem A' of \cite{JRW03} that $S_\I f$'s belong to $L^2(\Z)\subset L^\infty(\Z)\subset BMO_d(\Z)$ for any $f\in L^{\infty}_c(\Z)\subset L^2(\Z)$.

{\rm (ii)}. The case $\mathbb{T}$. Theorem \ref{main theorem 1} never be true in this case since $S_\I f$'s belong to $L^2(\mathbb{T})$, whence finit almost everywhere for any $f\in L^\infty(\mathbb{T})\subset L^2(\mathbb{T})$ by Theorem 2.2 in \cite{JKRW98}. An analog of Theorem \ref{main theorem 1} (ii) has been obtained in Remark 4.5 in \cite{JKRW98}.

{\rm (iii)}. The case $\R^n$. Let us explain in the case $n=2$ for simplifying the notations. An analog of Theorem \ref{main theorem 1} (i) is true by taking $f=\chi_{[0,\infty)\times[0,\infty)}$ and $Q_k(j,\ell)=[j-2^k,j)\times[\ell-2^k,\ell)$ using similar arguments. On the other hand, we can show an analog of Theorem \ref{main theorem 1} (ii) using similar calculations by considering $f$ and $\mathcal{Q}$ as follows. Let $f=\chi_{[1/2,1)}$. Let $Q_{\ell}=[1/2,1/2+1/2^{\ell})\times [1/2,1/2+1/2^{\ell})$. Let $P$ and $N$ denote two disjoint subsets of $[1/2,1)\times[1/2,1)$ such that for all $\ell>2$, $$|P\cap Q_{\ell}|=|N\cap Q_\ell|=2^{-2(\ell+1)}.$$
Take $Q_k(x,y)=(x-2^{k},x]\times (y-2^{k},y]$ for each $(x,y)\in N$, and  $Q_k(x,y)=(x,x+2^{k}]\times (y,y+2^{k}]$ for each $(x,y)\in P$.

 \end{remark}

\section{Proof of Theorem \ref{main theorem 2}}
This section is devoted to the proof of Theorem \ref{main theorem 2}. Moreover, we can verify that  the following argument work also in the case $\Z$, $\mathbb{T}$ and $\mathbb{R}^n$ with $n\geq2$. We leave the detailes for the interested readers.
\begin{proof}
It suffices to prove that there exist a positive constant $C$ such that
\begin{align}\label{l8cbmo}
\|Sf\|_{BMO_d}\leq C\|f\|_{\infty},\;\forall f\in L^{\infty}_c(\mathbb{R}).
\end{align}

We shall use the equivalent definition of $BMO_d$ norm.
\begin{align*}
\|g\|_{BMO_d}\simeq\sup_{I\;{\mathrm{dyadic}}}\inf_{a_{I}}\frac{1}{|I|}\int_{I}|g-a_{I}|.
\end{align*}

Give $f\in L^{\infty}(\R)$, and a dyadic interval $I$. We decompose $f$ as
 $f=f\un_{ I^* }+f\un_{\R\setminus I^*}=f_1+f_2$,
where $I^*$ is the cube with the same center as $I$ but three times the side length. We shall take $a_I=Sf_2(c_I)$ where $c_I$ is the center of $I$. Write $Sf-a_I$ as
$$Sf-a_I=Sf-Sf_2+Sf_2-a_I,$$
by triangle inequalities,
\begin{align*}
&\frac{1}{|I|}\int_{I}|Sf -a_{I}|\leq\frac{1}{|I|}\int_{I}|Sf-Sf_2| +\frac{1}{|I|}\int_I|Sf_2(x)- Sf_2(c_I)|dx\\
&\leq\frac{1}{|I|}\int_I\big(\sum_k\sup_{J\in\I_k}|M_{J+x}f_2(x)-\E_kf_2(x)-(M_{J+x}f_2(c_I)-\E_kf_2(c_I))|^{2}\big)^{\frac{1}{2}}\\
&+\frac{1}{|I|}\int_{I}|Sf_1|=(1)+(2).
\end{align*}
The first term $(1)$ is easily estimated by the fact that $S$ is of strong type $(2,2)$. Indeed,
\begin{align*}
(1)\leq\big(\frac{1}{|I|}\int_{I}|S(f_1)|^{2}\big)^{\frac{1}{2}}\leq
C\big(\frac{1}{|I|}\int_{I}|f_1|^{2}\big)^{\frac{1}{2}}\leq C\|f\|_{\infty}.
\end{align*}
The second term $(2)$ is controlled by a constant multiple of $\|f\|_{\infty}$ once we prove that for any $x\in I$ and any $J\in\mathcal{I}_k$,
$$\big(\sum_k|M_{J+x}f_2(x)-\E_kf_2(x)-(M_{J+x}f_2(c_I)-\E_kf_2(c_I))|^{2}\big)^{\frac{1}{2}}\leq C\|f\|_{\infty}.$$
If $2^k<|I|$, then $\E_kf_2$ is supported in $I$ and $J+x$ is contained in ${I}^*$ since
$$|x+y-c_I|\leq|x-c_I|+|y|\leq{1}/2|I|+2^{k+1}\leq 3|I|$$
for any $y\in J$. Then in this case, we get
$$M_{J+x}f_2(x)-\E_kf_2(x)-(M_{J+x}f_2(c_I)-\E_kf_2(c_I))=0,~\mbox{for any}~  x\in I.$$
Hence it suffices to consider the case $2^k\geq|I|$. Note that in this case, $I$ should be contained in some atom of $\sigma_k$, so $\E_kf_2(x)=\E_kf_2(c_I)$. On the other hand,
\begin{align*}
|M_{J+x}f_2(x)-&M_{J+x}f_2(c_I)|
=2^{-k}\big|\int_{J+x}f_2
-\int_{J+c_I}f_2\big|\\
&=2^{-k}\big|\int_{J+x\setminus J+c_I}f_2
-\int_{J+c_I\setminus J+x}f_2\big|\\
&\leq2^{-k}\int_{J+x\setminus J+c_I}|f_2|
+2^{-k}\int_{J+c_I\setminus J+x}|f_2|\\
&\leq2^{-k}|(J+x)\bigtriangleup(J+c_I)|\|f\|_{\infty}\leq C2^{-k}|I|\|f\|_{\infty}.
\end{align*}
The last inequality is due to the fact that $|(J+x)\bigtriangleup(J+c_I)|\leq C|x-c_I|\leq C|I|$.
Finally, the fact that $\ell^{2}$ norm is not bigger than $\ell^1$ norm implies
$$(2)\leq C|I|\|f\|_{\infty}\sum_{2^k\geq|I|}\frac{1}{2^{k}}\leq C\|f\|_{\infty}.$$
\end{proof}

\vskip30pt

\end{document}